\documentclass[12pt]{article}
\usepackage{amsxtra}
\usepackage{amssymb, amsmath}
\usepackage{amsthm}
\usepackage {hyperref}

\usepackage{amssymb}
\usepackage{amsmath}
\usepackage[dvips]{graphicx}

\begin{document}
\begin{center}
\textbf{\LARGE{Upside down Magic, Bimagic, Palindromic Squares and  Pythagoras Theorem on a Palindromic Day - 11.02.2011}}
\end{center}

\bigskip
\begin{center}
\textbf{\large{Inder Jeet Taneja}}\\
Departamento de Matem\'{a}tica\\
Universidade Federal de Santa Catarina\\
88.040-900 Florian\'{o}polis, SC, Brazil.\\
\textit{e-mail: taneja@mtm.ufsc.br\\
http://www.mtm.ufsc.br/$\sim $taneja}
\end{center}

\begin{abstract}
\textit{In this short note we have produced different kinds of upside down magic squares based on a palindromic day 11.02.2011. In this day appear only the algorisms 0, 1 and 2. Some of the magic squares are bimagic and some are palindromic. Magic sums of the magic squares of order 3$\times$3, 4$\times$4 and 5$\times$5 satisfies the Pythagoras theorem. Three different kinds of bimagic squares of order 9$\times$9 are also produced. The bimagic square of order 9$\times$9 with 8 digits is palindromic numbers. We have given bimagic squares of order 16$\times$16 and 25$\times$25, where the magic sum S1 in both the cases is same. In order to make these magic squares upside down, i.e., $180^{0}$  degree rotation,  we have used the numbers in the digital form. All these magic square are only with three digits, 0, 1 and 2 appearing in the day 11.02.2011.}

\end{abstract}

\section{Introduction}

It interesting to observe that the day 11.02.2011 is palindromic and has only three digits 0, 1 and 2. A similar kind of palindromic day shall also appear next year 21.02.2012 having the same three  digits. In this paper our interest is to produce upside down magic squares, bimagic squares and palindromic magic squares using only these three algorisms, 0, 1 and 2. A similar kind of study can be seen in another author's work on the day October 10, 2010 \cite{tan3}. Using these three digits we have \cite{tan6} made equivalence with classical magic squares, the one is \textit{"Lo-Shu"}  magic square of order 3$\times$3 and another is \textit{"Khajurao"} magic squares of order 4$\times$4.

\bigskip
Before we proceed, here below are some basic definitions:
\begin{itemize}
\item [(i)]  A \textbf{magic square} is a collection of numbers put as a square matrix, where the sum of element of each row, sum of element of each column and sum of each element of two principal diagonals have the same sum.  For simplicity, let us write it as \textbf{S1}.

\item  [(ii)]\textbf{Bimagic square} is a magic square where the sum of square of each element of rows, columns and two principal diagonals are the same. For simplicity, let us write it as \textbf{S2}.

\item [(iii)] \textbf{Upside down}, i.e., if we rotate it to $180^{0}$ degree it remains the same.

\item [(iv)] \textbf{Mirror looking}, i.e., if we put it in front of mirror or see from the other side of the glass, or see on the other side of the paper, it always remains the magic square.

\item [(v)] \textbf{Universal magic squares}, i.e., magic squares having the property of upside down and mirror looking are considered \textit{universal magic squares}.
\end{itemize}

In this short note we have produced different kinds of upside down magic squares using only the algorisms 0, 1 and 2. Some of the magic squares are bimagic and some are palindromic. Magic sum of the magic squares of order 3$\times$3, 4$\times$4 and 5$\times$5 satisfies the Pythagoras theorem. Bimagic squares of order 9x9 are produced with 4, 6 and 8 digits. The bimagic square of order 9$\times$9 with 8 digits is palindromic while with 6 digits is a combination of palindromic numbers. We have given bimagic squares of order 16$\times$16 and 25$\times$25, where the magic sum S1 is the same. In order to make these magic squares upside down we have used the numbers in the digital form.

\bigskip
All these magic square are only with three digits, 0, 1 and 2. In order to do so, we have used the numbers in the digital form:

\begin{center}
\includegraphics[bb=0mm 0mm 208mm 296mm, width=20.8mm, height=5.0mm, viewport=3mm 4mm 205mm 292mm]{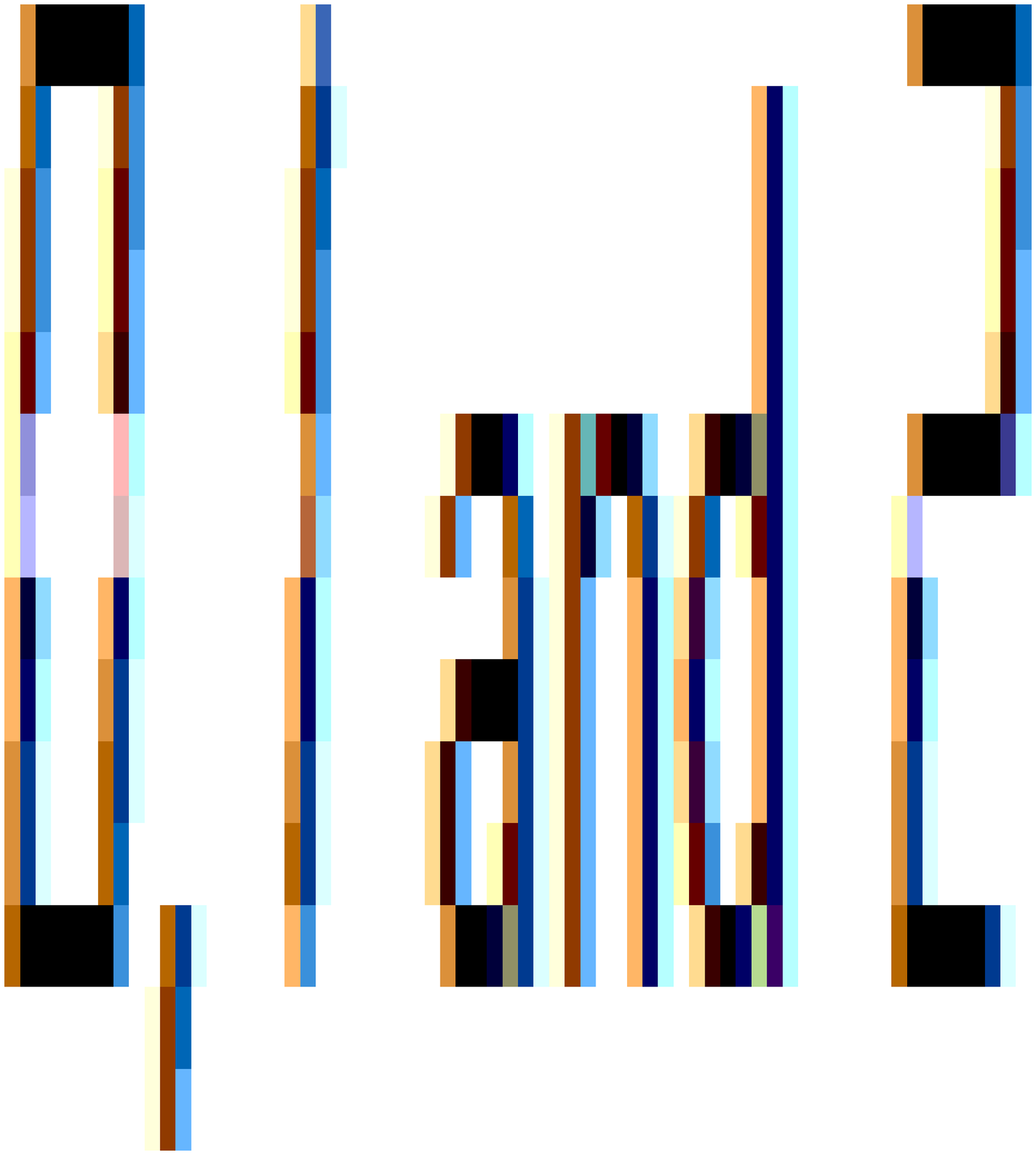}
\end{center}

These digits generally appear in watches, elevators, etc.  We observe that the above three digits are rotatable to $180^{0}$, and remains the same.  These also be considered as universal, because in the mirror 2 becomes 5, while 0 and 1 remains the same. In these situations the magic sums are different. This we leave to reader to verify.

\section{Upside Down Magic Squares and the Pythagoras Theorem}

In this section we shall present magic squares of order 3$\times$3, 4$\times$4 and 5$\times$5 having only the three digits 0, 1 and 2 in the digital form. Interesting the magic sum S1, in this case satisfies the Pythagoras theorem.

\bigskip
\noindent
\textbf{$\bullet$ Magic squares of order 3$\times$3}

\bigskip
Here below are two magic squares of order 3$\times$3 with $S1_{3\times 3} :=33$ and $S1_{3\times 3} :=3333$ respectively. The first one is with two digits combinations while the second one is with four digits combinations:

\begin{center}
{\includegraphics[bb=0mm 0mm 208mm 296mm, width=24.6mm, height=18.0mm, viewport=3mm 4mm 205mm 292mm]{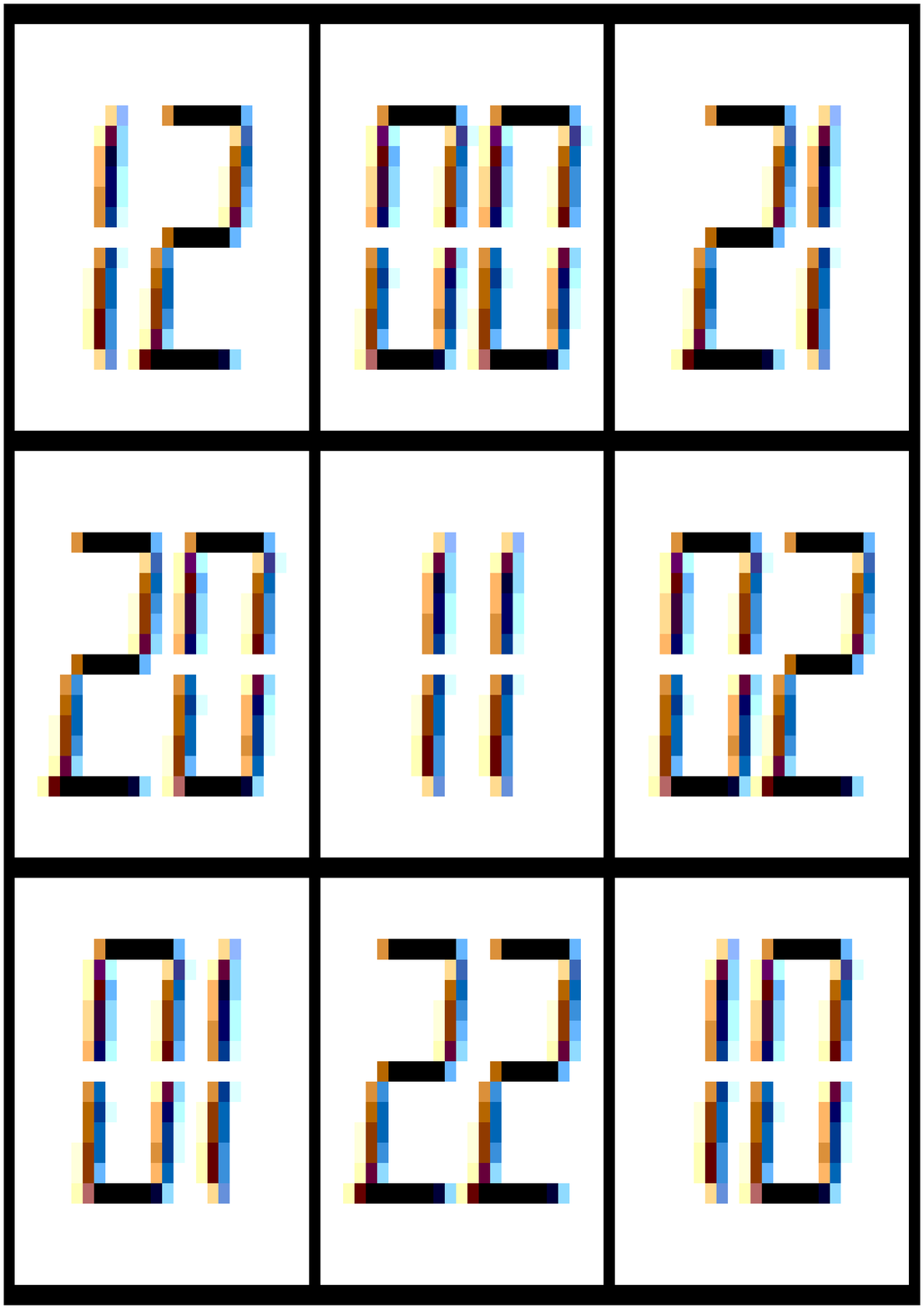}}
\end{center}

\begin{center}
{\includegraphics[bb=0mm 0mm 208mm 296mm, width=36.1mm, height=24.4mm, viewport=3mm 4mm 205mm 292mm]{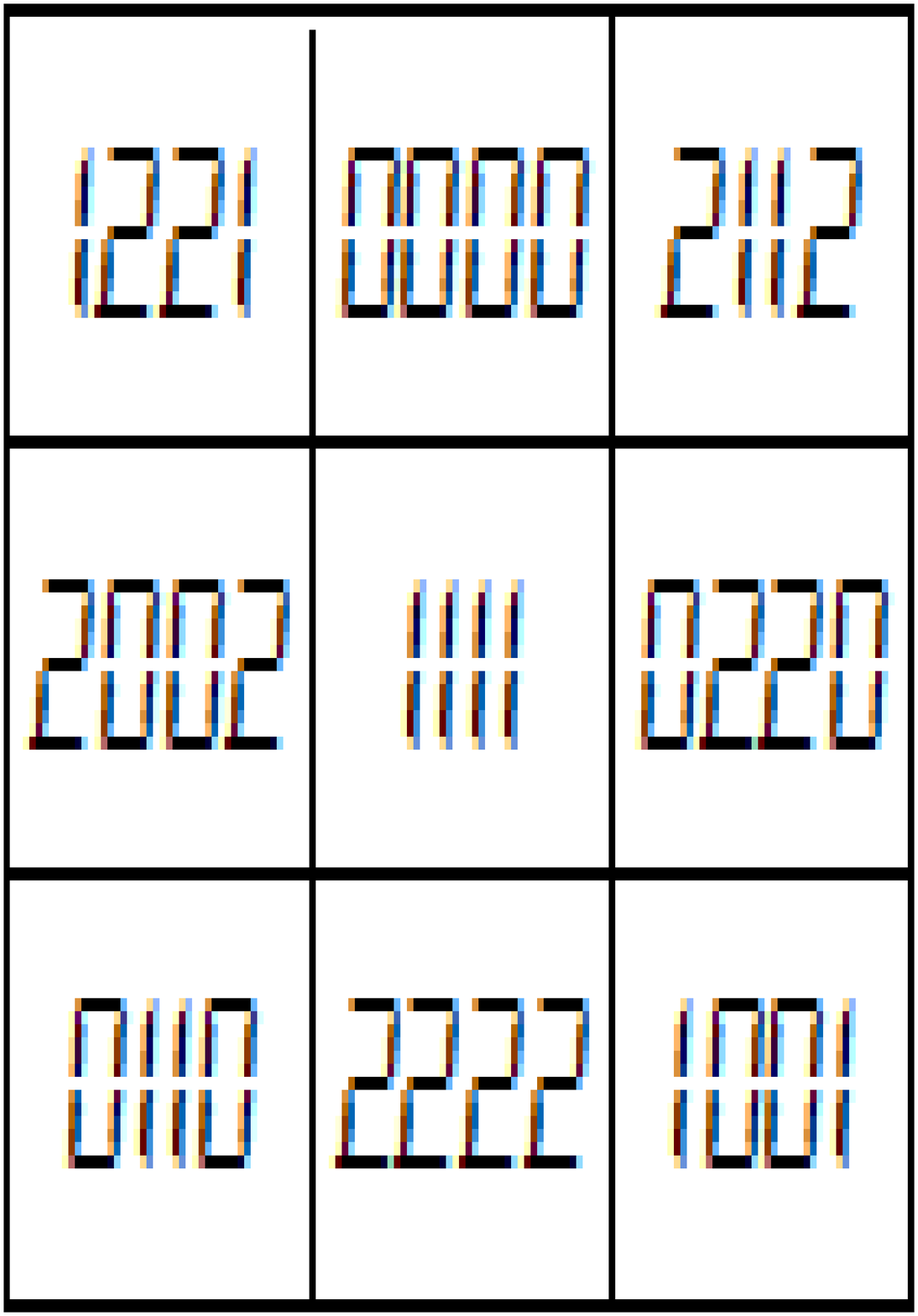}}
\end{center}

We observe from the second magic square that it is \textbf{palindromic}. In order to have upside down we have considered 110, 220 as 0110, 0220 to be symmetry in the result.

\bigskip
\noindent
\textbf{$\bullet$ Magic squares of order 4$\times$4}

\bigskip
Here below is a magic square of order 4x4 with $S1_{4\times 4} :=4444$

\begin{center}
\includegraphics[bb=0mm 0mm 208mm 296mm, width=49.6mm, height=37.8mm, viewport=3mm 4mm 205mm 292mm]{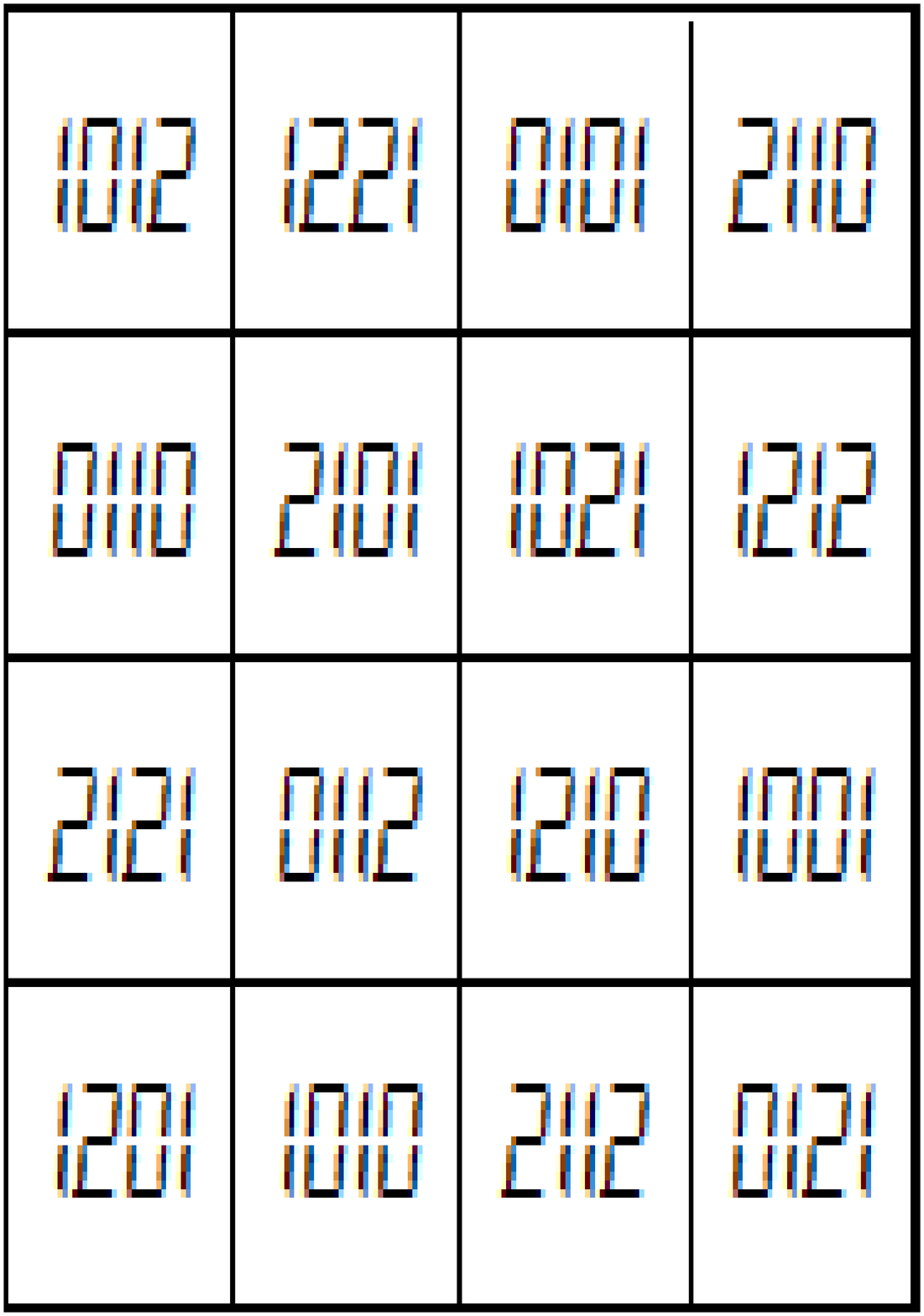}
\end{center}

\bigskip
\noindent
\textbf{$\bullet$ Magic squares of order 5$\times$5}

\bigskip
Here below is a magic square of order 5$\times$5 with $S1_{5\times 5} :=5555$

\begin{center}
\includegraphics[bb=0mm 0mm 208mm 296mm, width=65.3mm, height=49.6mm, viewport=3mm 4mm 205mm 292mm]{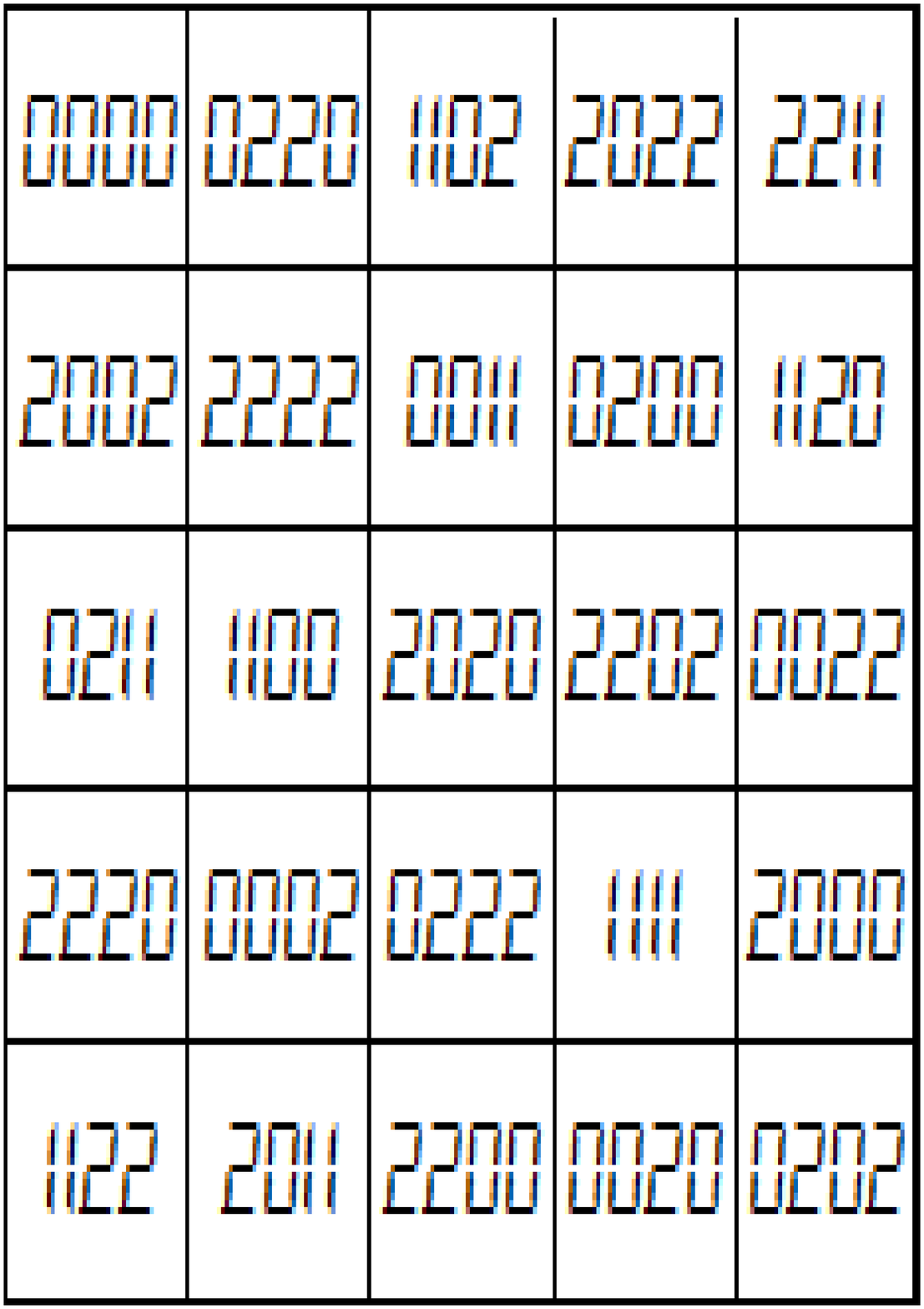}
\end{center}

\bigskip
The above magic square is pan diagonal.

\subsection{Pythagoras Theorem}

From the above magic squares of orders 3$\times$3, 4$\times$4 and 5$\times$5 with four digits, we have the following result:

\[\left(S1_{3\times 3} \right)^{2} +\left(S1_{4\times 4} \right)^{2} =\left(S1_{5\times 5} \right)^{2} ,\]

i.e.,

\[3333^{2} +4444^{2} =5555^{2} \]

i.e.,

\[11108889+19749136=30858025.\]

This means that if we consider square of any line, or column or principle diagonal, from one of the above magic squares we shall always have the same value, for example,

\[\begin{array}{l} {\left(1221+1111+1001\right)^{2} +\left(1012+2101+1210+0121\right)^{2} } \\ {\quad \quad =\left(2002+2222+0011+0200+1120\right)^{2} {\kern 1pt} } \end{array}\]

The above result gives us an upside down \textbf{Pythagoric equation} when we use the digital letters:

\begin{center}
\includegraphics[bb=0mm 0mm 208mm 296mm, width=123.1mm, height=10.5mm, viewport=3mm 4mm 205mm 292mm]{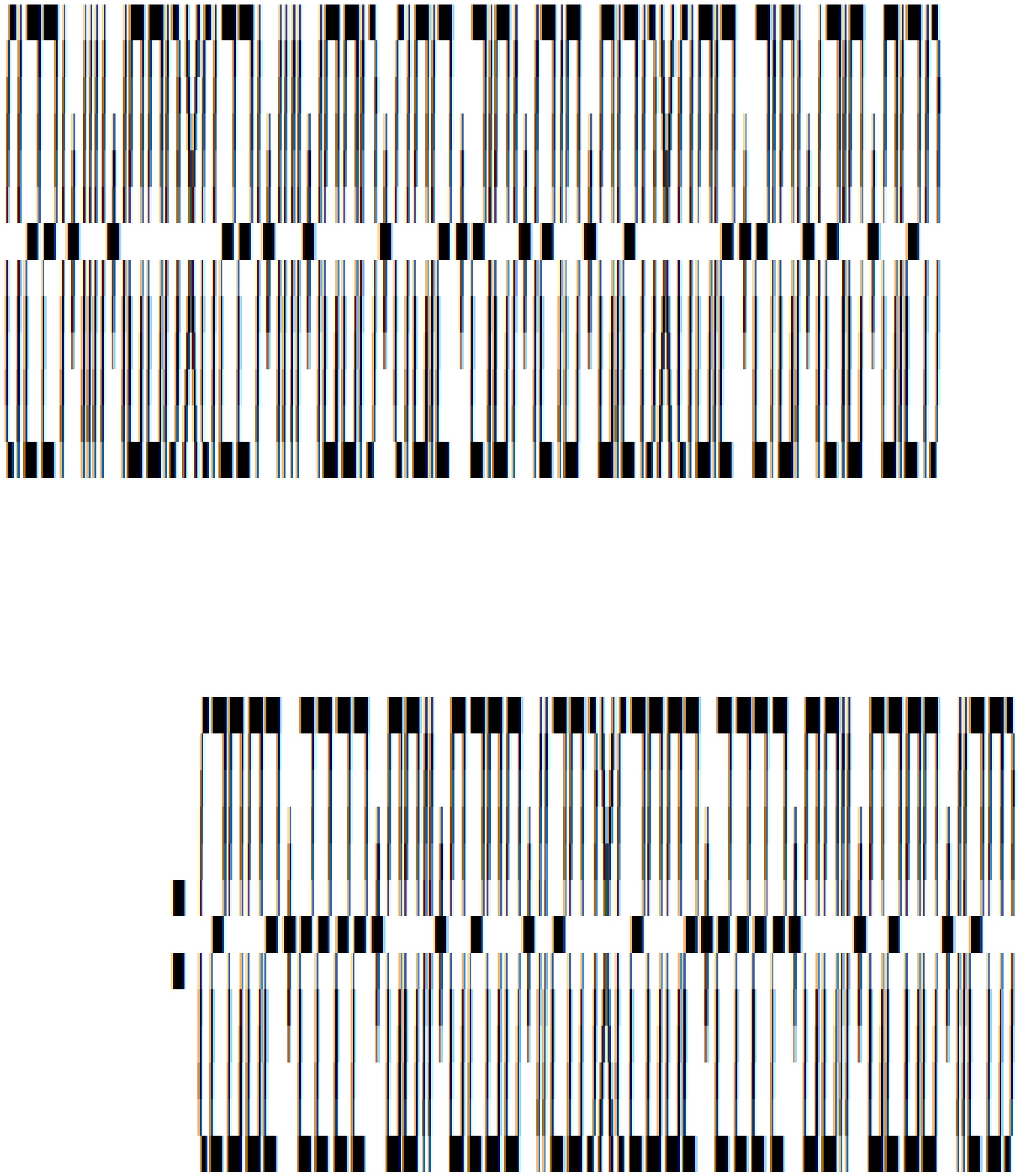}
\end{center}

$180^{0}$ degrees rotation:

\begin{center}
\includegraphics[bb=0mm 0mm 208mm 296mm, width=123.1mm, height=10.5mm, viewport=3mm 4mm 205mm 292mm]{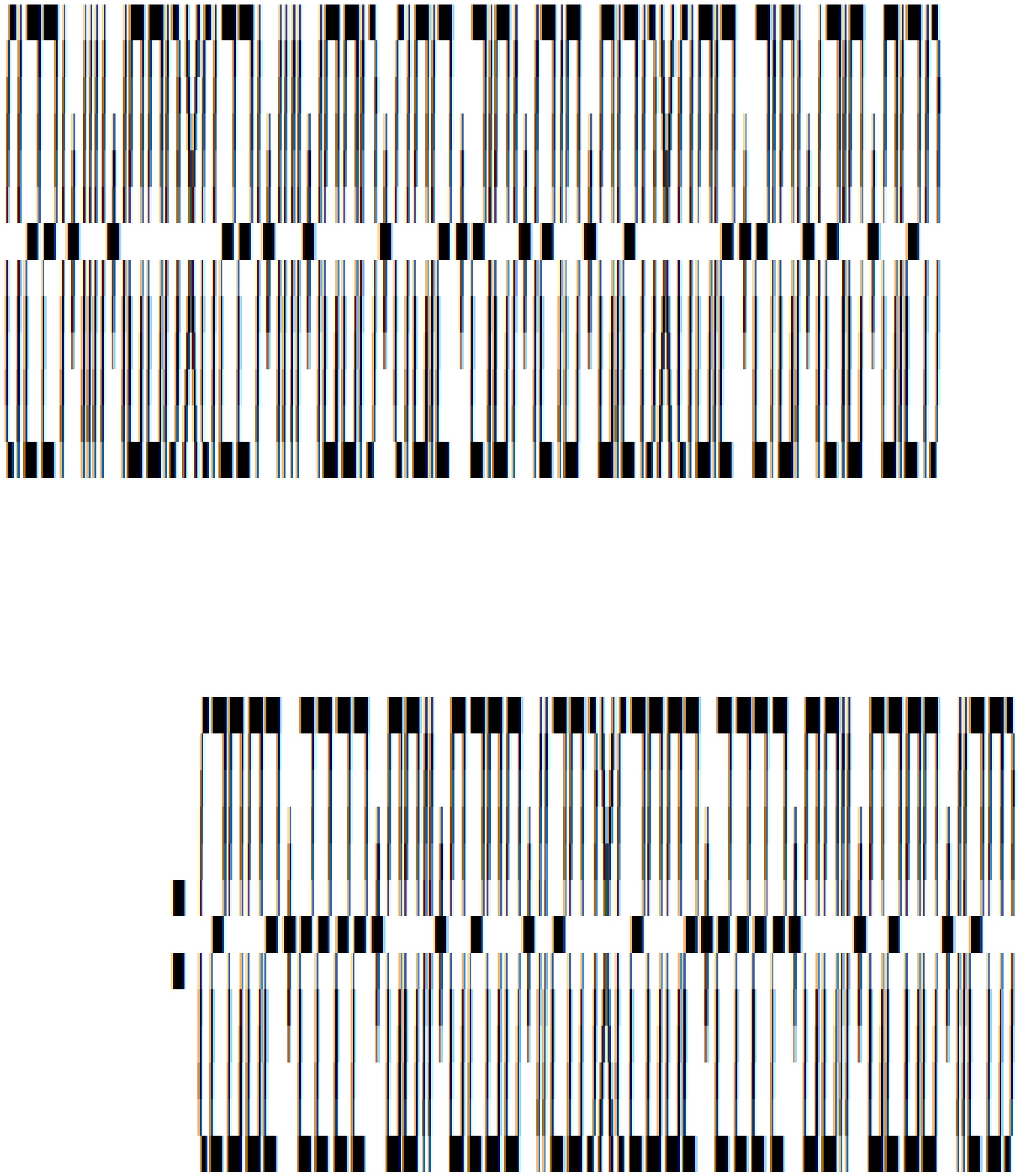}
\end{center}

We observe from the second equation that the numbers are different, but the sum is same.

\section{Upside down bimagic squares of order 9x9}

\bigskip
Here below are three bimagic squares of order 9x9. The first one is with four digits combinations, the second one is with six digits combinations and the third one is with eight digits combinations.

\bigskip
\noindent
\textbf{$\bullet$ Bimagic squares of order 9x9 with four digits}

\begin{center}
\includegraphics[bb=0mm 0mm 208mm 296mm, width=108.7mm, height=59.8mm, viewport=3mm 4mm 205mm 292mm]{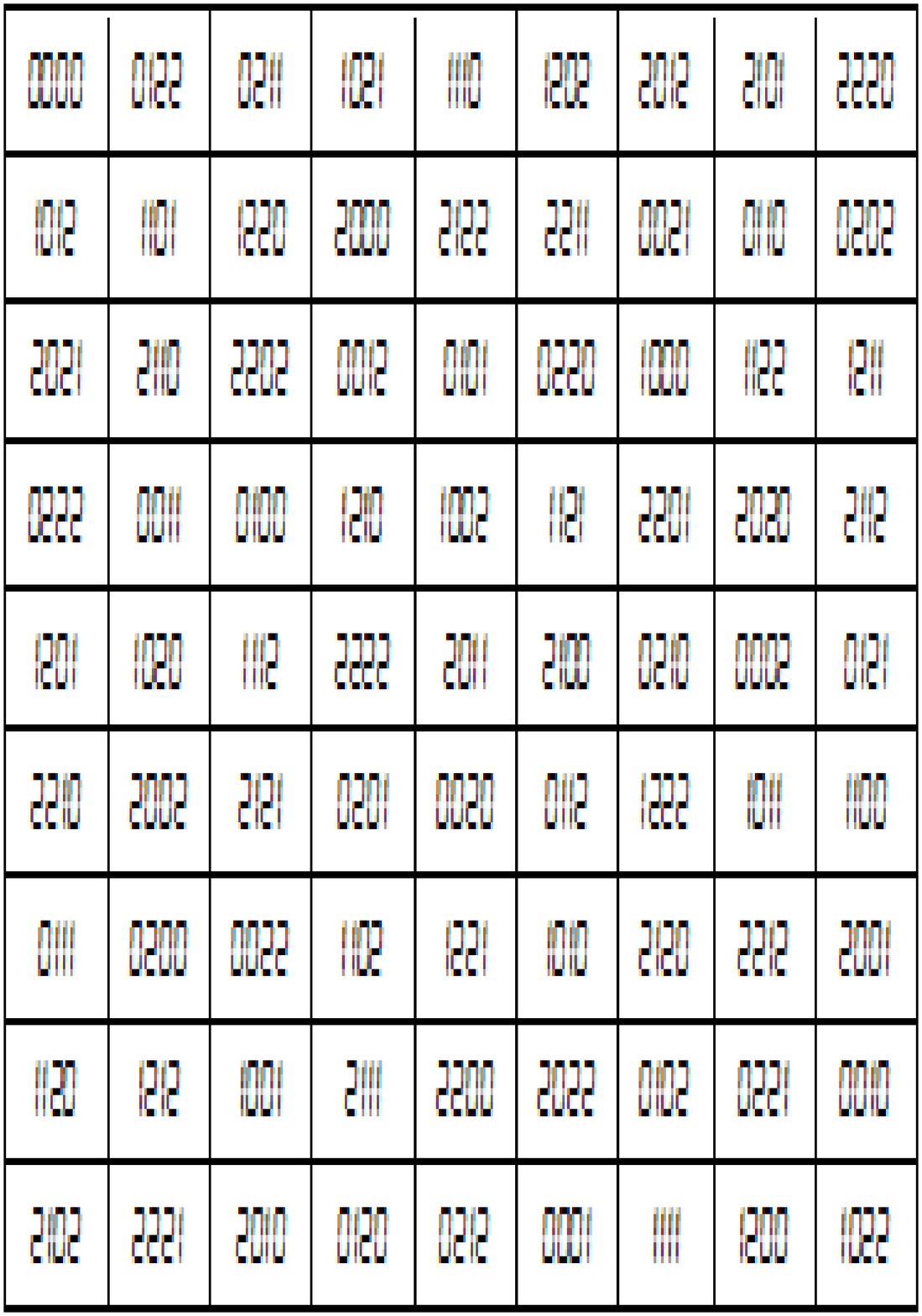}
\end{center}

\noindent
$S1_{9\times 9} :=9999$\\
$S2_{9\times 9} :={\rm 17169395}$

\bigskip
Also we have sum of each block of order 3x3 is 9999 and the squares of sum of each term in each block of order 3x3 is also 17169495. We observe that the above bimagic square still has the day's number 11.02.2011 in two parts 1102 and 2011.

\bigskip
\noindent
\textbf{$\bullet$ Bimagic squares of order 9$\times$9 with six digits}

\begin{center}
\includegraphics[bb=0mm 0mm 208mm 296mm, width=129.2mm, height=60.0mm, viewport=3mm 4mm 205mm 292mm]{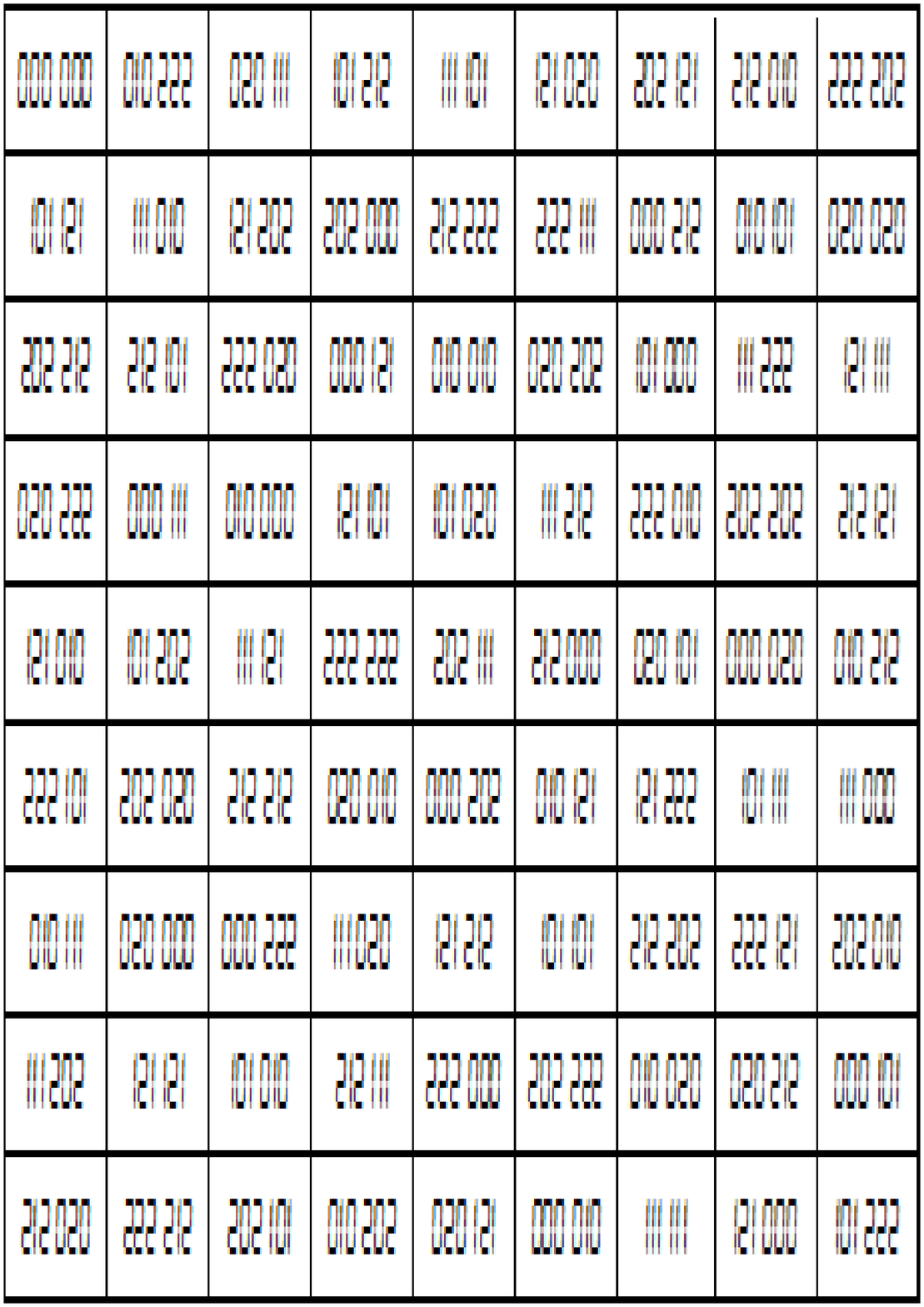}
\end{center}

\bigskip
\noindent
$S1_{9\times 9} :=999999$\\
$S2_{9\times 9} :={\rm 17291695}0{\rm 695}$,

\bigskip
Also we have sum of each block of order 3x3 is 999999 and square of sum of each term in each block of order 3x3 is also 172916950695. Here we observe that each number is a composition of three digit palindromic numbers.

\bigskip
\noindent
 \textbf{$\bullet$ Bimagic squares of order 9x9 with eight digits}

\begin{center}
\includegraphics[bb=0mm 0mm 208mm 296mm, width=143.7mm, height=60.0mm, viewport=3mm 4mm 205mm 292mm]{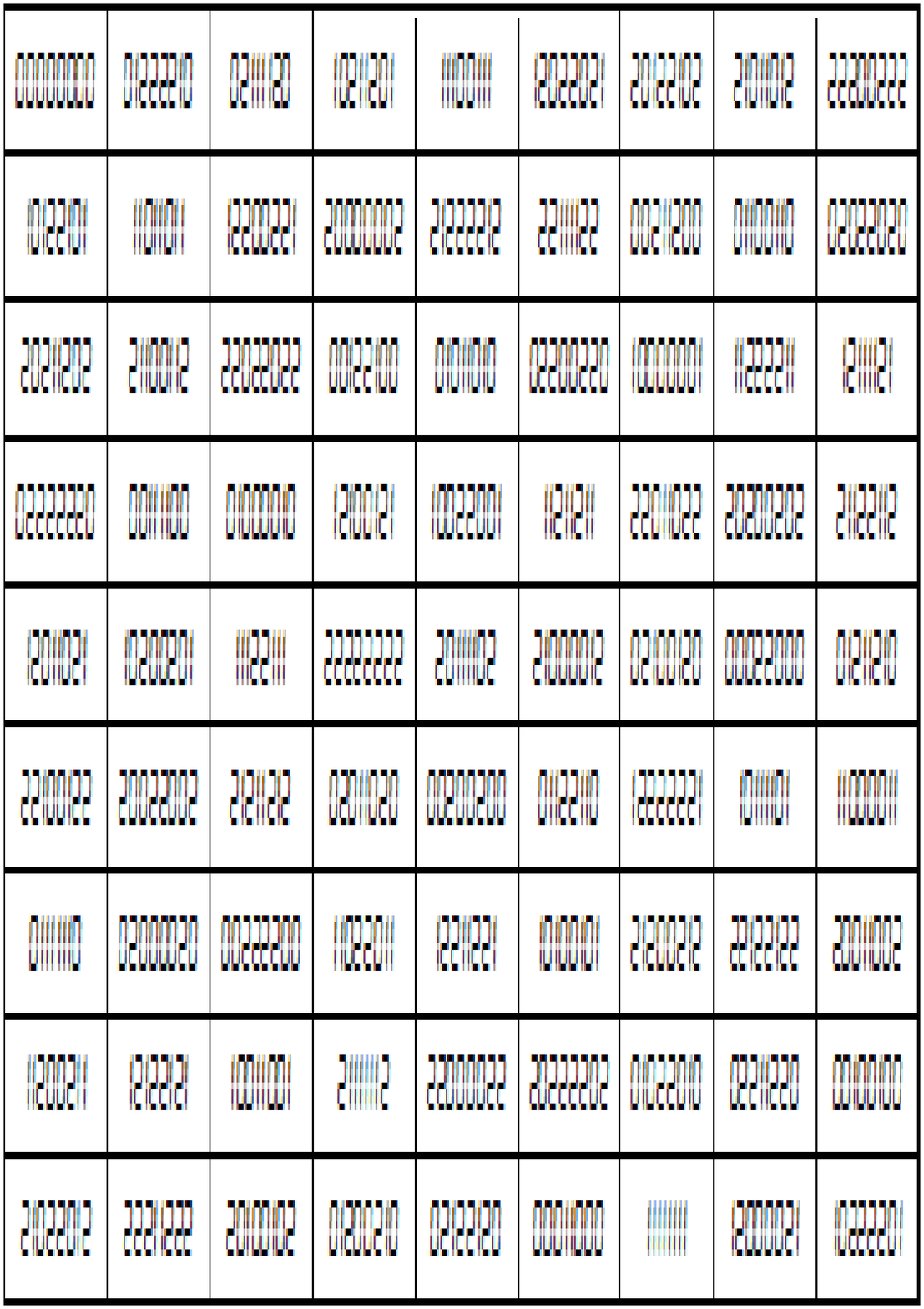}
\end{center}

\noindent
$S1_{9\times 9} :=99999999$\\
$S2_{9\times 9} :={\rm 1717172174949490}$

\bigskip
Also we have sum of each block of order 3$\times$3 is 99999999 and square of sum of each term in each block of order 3$\times$3 is 1717172174949490. Here the numbers are palindromic including the day's number:  \textbf{11.02.2011}.

\section{Upside down bimagic squares of order 16$\times$16 and 25$\times$25 with same magic sum}

\bigskip
Here below are bimagic squares of order 16$\times$16 and 25$\times$25. Both these magic square have the same magic sum, $S1_{16\times 16} =S1_{25\times 25} =222222220$.

\bigskip
\noindent
\textbf{$\bullet$ Bimagic square of order 16$\times$16 }

\begin{center}
\includegraphics[bb=0mm 0mm 208mm 296mm, width=159.7mm, height=76.0mm, viewport=3mm 4mm 205mm 292mm]{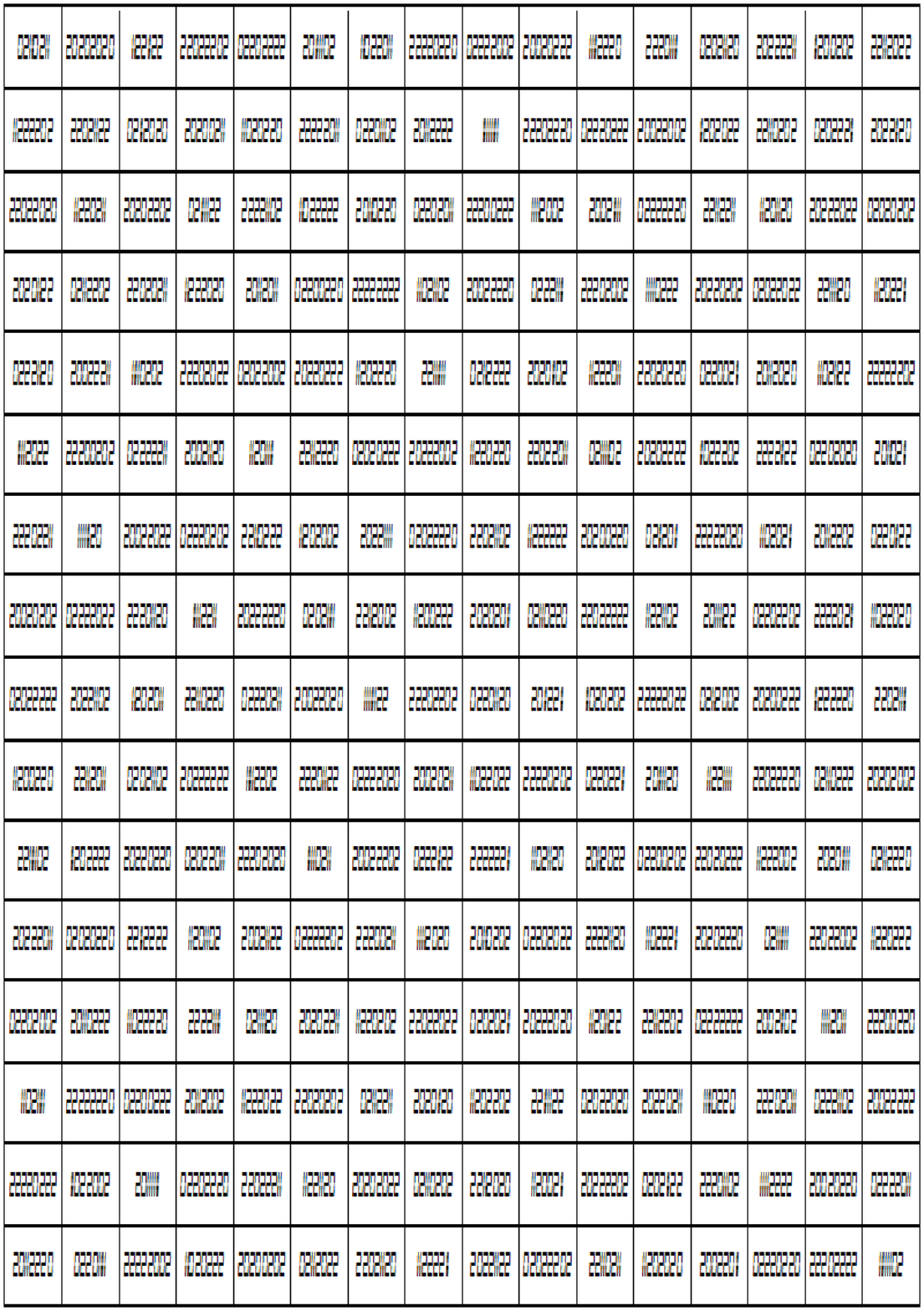}
\end{center}

\noindent
$S1_{16\times 16} :=222222220$\\
$S2_{16\times 16} :={\rm 4}0{\rm 9752}0{\rm 8}0{\rm 1469}0{\rm 4}0$.

\bigskip
Here each block of order 4$\times$4 is a magic square of sum$S1_{4\times 4} :=55555555$. Also we observe that the numbers are not palindromic, but still include the day's number:  \textbf{11.02.2011. }

\bigskip
\noindent
\textbf{$\bullet$ Pan diagonal bimagic square of order 25x25 }

\begin{center}
\includegraphics[bb=0mm 0mm 208mm 296mm, width=78.0mm, height=18.1mm, viewport=3mm 4mm 205mm 292mm]{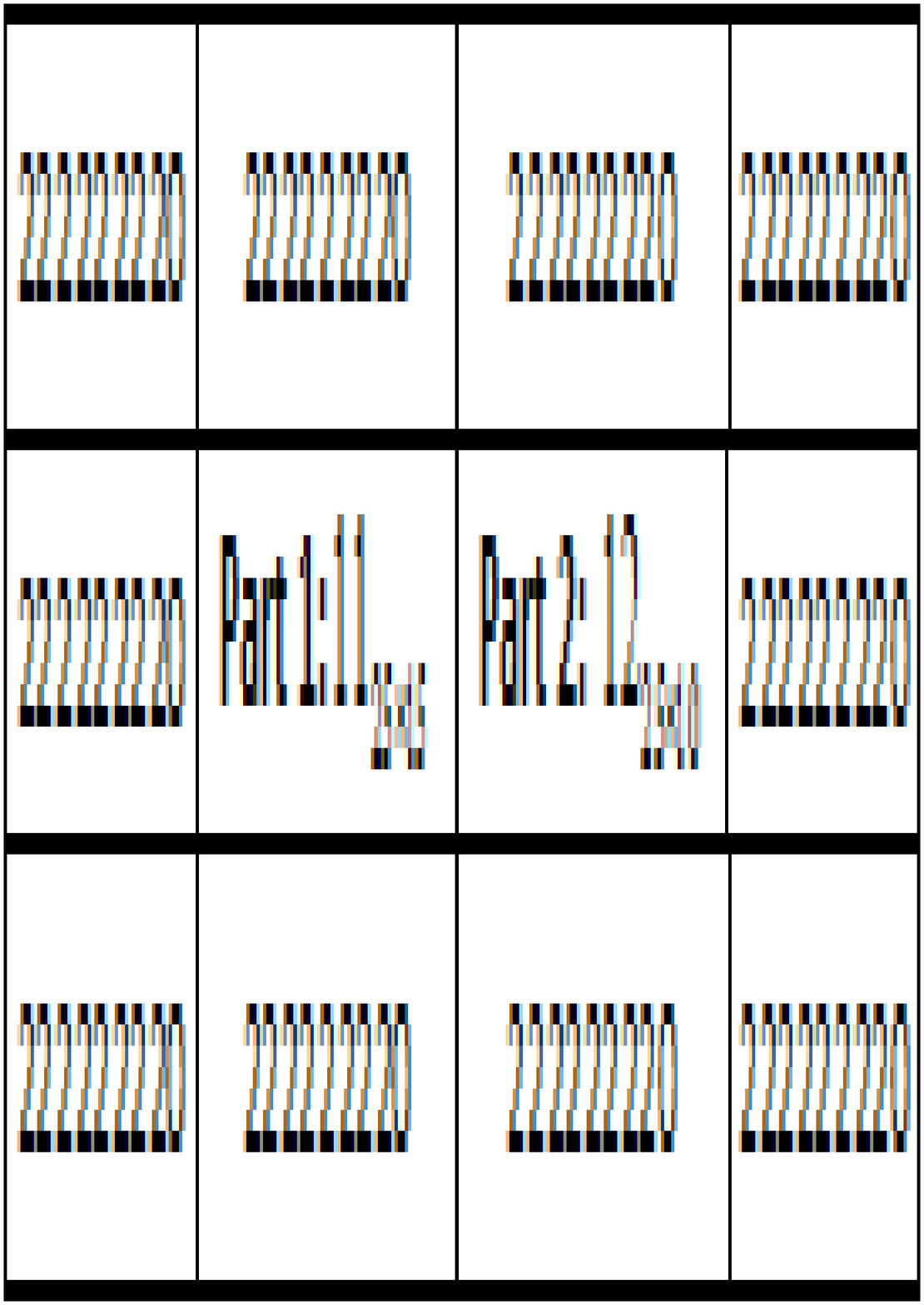}
\end{center}

\bigskip
\noindent
\textbf{Part 1:$11_{25\times 15} $}

\begin{center}
\includegraphics[bb=0mm 0mm 208mm 296mm, width=160.0mm, height=111.3mm, viewport=3mm 4mm 205mm 292mm]{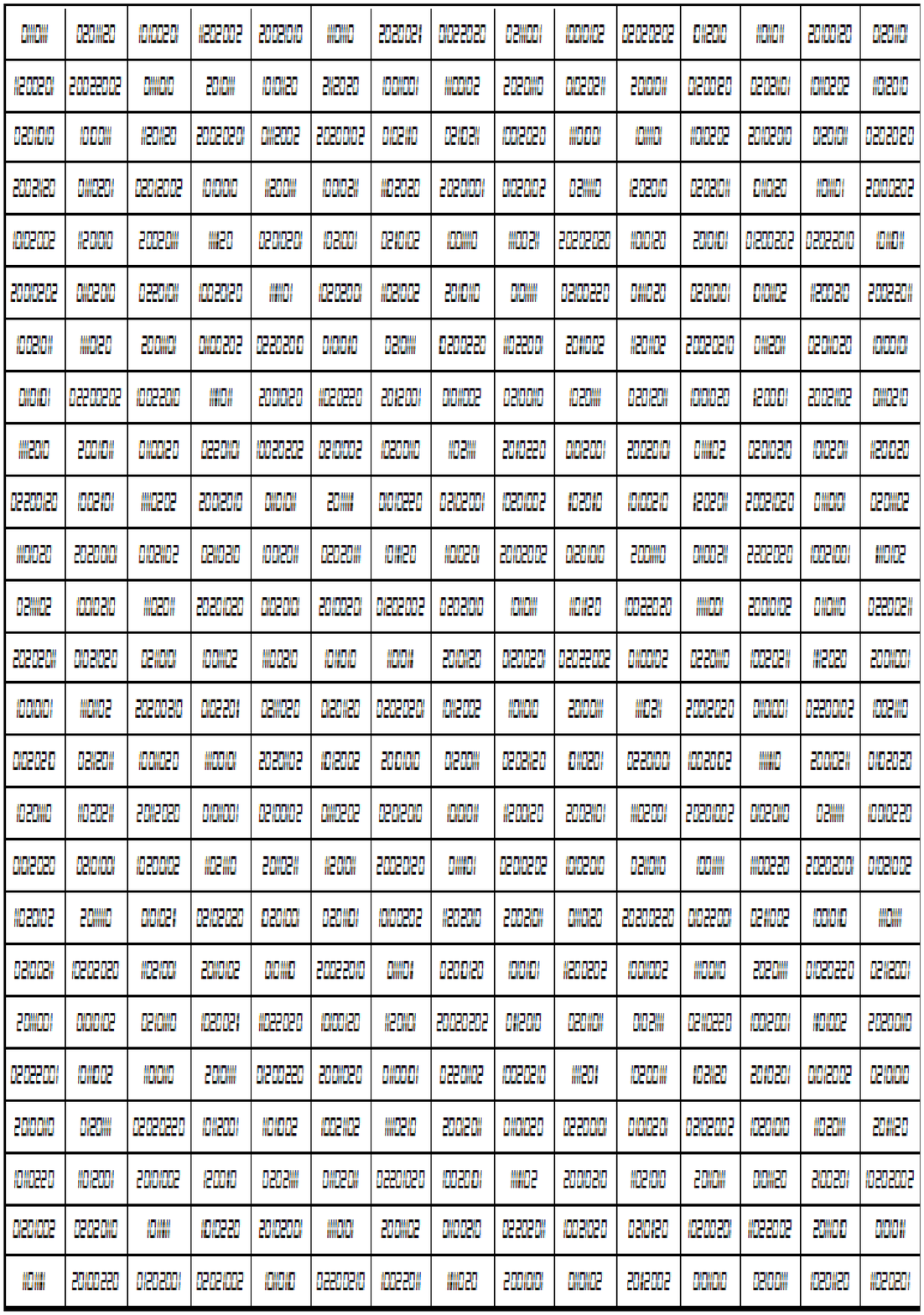}
\end{center}

\bigskip
\noindent
\textbf{Part 2: $12_{25\times 10} $}

\begin{center}
\includegraphics[bb=0mm 0mm 208mm 296mm, width=121.7mm, height=111.4mm, viewport=3mm 4mm 205mm 292mm]{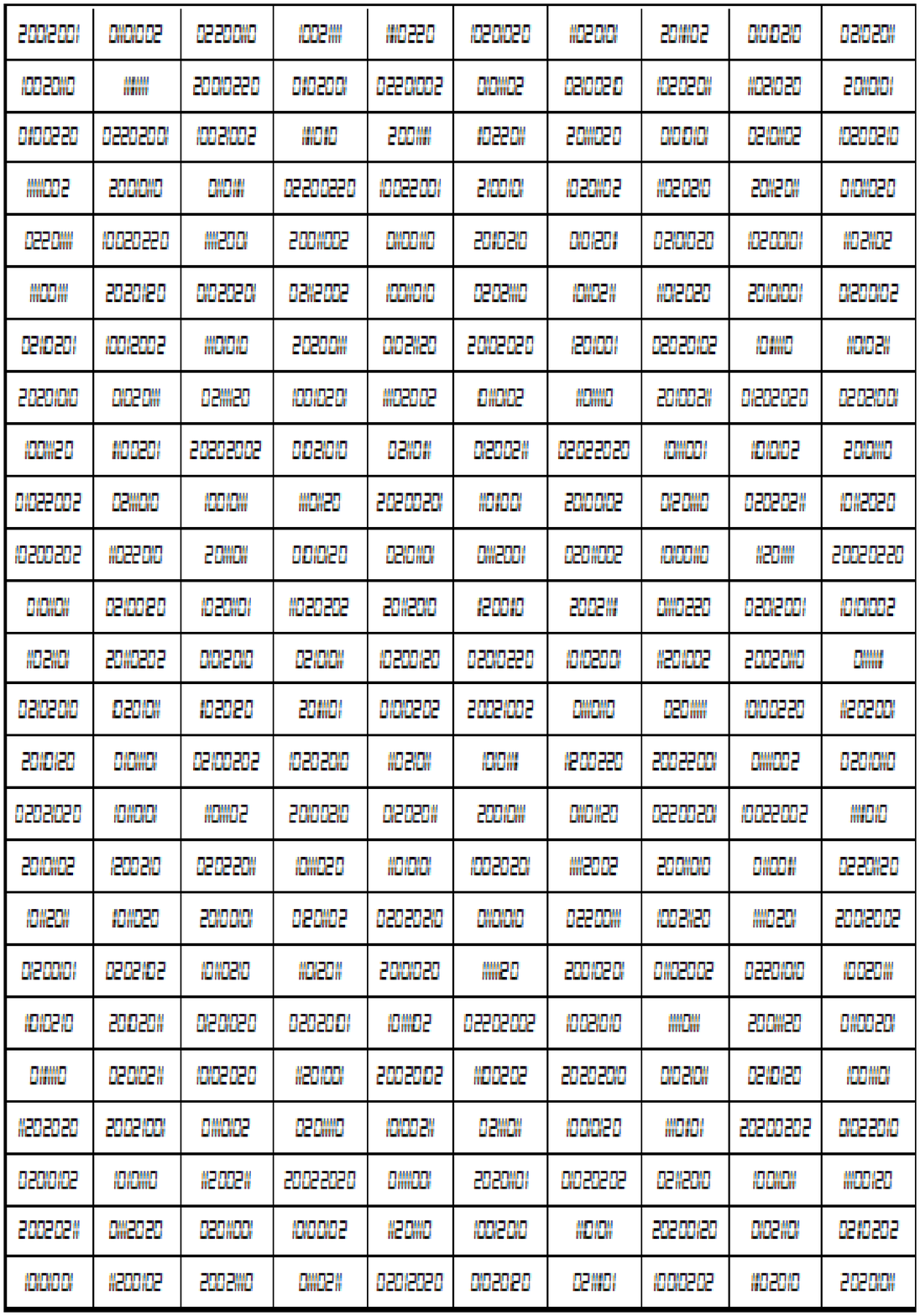}
\end{center}

\noindent
$S1_{25\times 25} :=222222220$\\
$S2_{25\times 25} :=3169428014410330$

\bigskip
Here each block of order 5$\times$5 is a magic square of sum  $S1_{5\times 5} :=44444444$. We observe that the numbers are not palindromic, but still include the day's number:  \textbf{11.02.2011. }

\section{Final Comments}
In this short paper we have brought magic squares of different kinds using only the digits 0, 1 and 2 appearing in a palindromic day - 11.02.2011. Another palindromic day having the same digits shall also appear next year - 21.02.2012. The magic squares obtained are upside down, i.e., when we make a rotation of $180^{0}$ degrees they still remains the magic squares. This happens using the letters in the digital. Some of the magic squares are palindromic. Bimagic squares are of order 9$\times$9, 16$\times$16 and 25$\times$25. The numbers 9, 16 and 25 remember us a Pythagoras theorem, i.e.,  $3^{2} + 4^{2} = 5^{2}$. We have produced magic squares of orders 3$\times$3, 4$\times$4 and 5$\times$5 using only three digits 0, 1 and 2 and the sum S1 satisfies the Pythagoras theorem, i.e.,  $\left(S1_{3\times 3} \right)^{2} +\left(S1_{4\times 4} \right)^{2} =\left(S1_{5\times 5} \right)^{2}$.  Interestingly, the sum S1 in case of bimagic squares of order 16$\times$16 and 25$\times$25 is the same and can also be made upside down. Most of the magic squares have the palindromic number 11.02.2011. The digits 0, 1 and 2 also appears in many others days during the years 2010, 2011 and 2012. In another work \cite{tan6}, we have brought equivalent versions of two classical magic squares of order 3$\times$3 and 4$\times$4 using only these three digits 0, 1 and 2. For more studies on magic squares see the references below.

\begin{center}
---------------------------
\end{center}

\end{document}